\documentclass{article}
\usepackage[utf8]{inputenc}
\usepackage[english]{babel}
\usepackage{graphicx}
\usepackage{amsmath, accents}
\usepackage{amssymb}
\usepackage{amsfonts}
\usepackage{amsthm}
\usepackage{mathrsfs}
\usepackage{latexsym}
\usepackage{subfigure}
\usepackage[ruled,vlined]{algorithm2e}
\usepackage{color}
\usepackage{xcolor}
\usepackage{enumitem}
\usepackage{hyperref,authblk}

\theoremstyle{definition}
\newtheorem{teo}{Theorem}

\newtheorem{rem}[teo]{Remark}

\begin{document}

\title{A ``right" path to cyclic polygons}

\author[1]{Paolo Dulio\thanks{paolo.dulio@polimi.it}\footnote{Corresponding author}}

\author[1]{Enrico Laeng\thanks{enrico.laeng@polimi.it}}

\affil[1]{Dipartimento di Matematica ``F.~Brioschi'', Politecnico di Milano, Piazza Leonardo da Vinci $32$, I-$20133$ Milano, Italy}

\date{ }


%

\thispagestyle{empty}

\maketitle

\begin{abstract}
It is well known that Heron's theorem provides an explicit formula for the area of a triangle, as a symmetric function of the lengths of its sides. It has been extended by Brahmagupta to quadrilaterals inscribed in a circle (cyclic quadrilaterals). A natural problem is trying to further generalize the result to cyclic polygons with a larger number of edges, which, surprisingly, has revealed to be far from simple. In this paper we investigate such a problem by following a new and elementary approach. We start from the simple observation that the incircle of a right triangle touches its hypothenuse in a point that splits it into two segments, the product of whose lengths equals the area of the triangle. From this curious fact we derive in a few lines: an unusual proof of the Pythagoras' theorem, Heron's theorem for right triangles,  Heron's theorem for general triangles, and Brahmagupta's theorem for cyclic quadrangles. This suggests that cutting the edges of a cyclic polygon by means of suitable points should be the ``right" working method. Indeed, following this idea, we obtain an explicit formula for the area of any convex cyclic polygon, as a symmetric function of the segments split on its edges by the incircles of a triangulation. We also show that such a symmetry can be rediscovered in Heron’s and Brahmagupta’s results, which consequently represent special cases of the general provided formula.\smallskip

\noindent MSC: 52A10;52A38

\end{abstract}\bigskip

\textbf{Keywords:} Area; cyclic polygon; incircle; inradius.


\section{Introduction}

A natural and largely considered question in convex geometry is the determination of the area $A$ of a convex polygon as a function of the lengths of its sides. The problem goes back to Heron of Alexandria, that was able to solve the problem in the case of a triangle. Later, in the seventh century, Brahmagupta extended the result to cyclic quadrilaterals, namely to quadrilaterals inscribed in a circle (see for instance \cite{CG}). Several results concerning the geometry of cyclic polygons have been obtained in different areas of research (see \cite{CK,D,DG,KSS,WX}), which points out a general interest for such geometric objects. It is therefore natural trying to further extend to cyclic polygons with a larger number of edges the nice and ancient formulae by Heron and Brahmagupta. Surprisingly, this has revealed to be far from simple. In \cite{Robbins1,Robbins2}  an algebraic formulation of the problem led D.P. Robbins to find formulae for cyclic pentagons and cyclic hexagons. It was observed that Heron and Brahmagupta's formulae can be restated in a form where $16A^2$ represents a monic polynomial whose coefficients are symmetric polynomials in the squares of the sides. This generalizes to cyclic pentagons and hexagons, where the polynomial have degree $7$, but the formulae, even if  holding also in the non convex case, do not provide an easy explicit form for the area (see also \cite{PAK} for interesting comments and remarks). Formulae of the same kind have been conjectured \cite{Robbins1},  and later proved \cite{MRR}, even for heptagons and octagons, also illuminating some mysterious features of Robbins' formulas for the areas of cyclic pentagons and hexagons (see also \cite{V} for further detail on Robbin's conjectures). The resulting formulae are interesting, but are very complex and do not seem to provide a general picture that could be easily generalizable to polygons with an arbitrary large number of edges.

In this article we follow a different approach, which leads to a complete solution of the considered problem.  The leitmotif of our paper is to point out that the role played by the edges in Heron’s and Brahmagupta’s theorems must be replaced by the segments cut on the edges of a polygon by the incircles of the triangles of a triangulation of the polygon.
First of all, we show that Pytagoras', Heron's and Brahmagupta's theorems can be linked together thanks to a simple result concerning the area of a right triangle. Then, we generalize Heron's and Brahmagupta's results to a symmetric coordinate free formula that holds true for any cyclic polygon.

\section{Heron's formula}

\begin{center}
\begin{figure}[htb]
\centering
\includegraphics[scale=0.8,bb=118 459 488 632,clip=true]{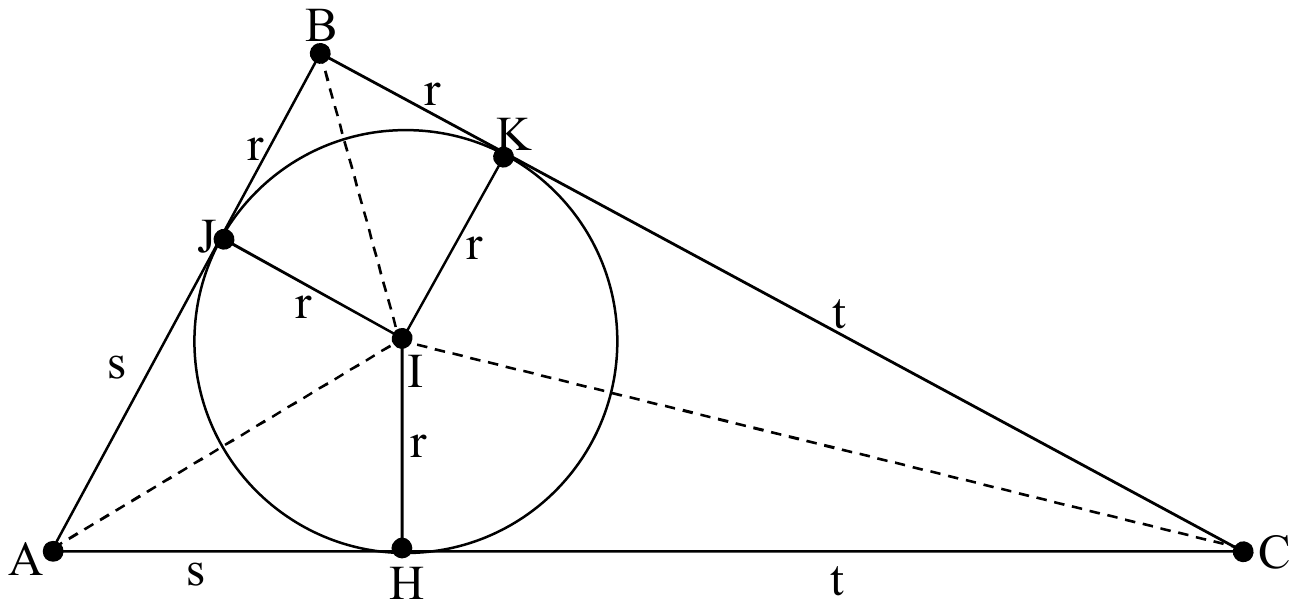}\caption{A right triangle $ABC$.} \label{fig:Erone}
\end{figure}
\end{center}

Let $ABC$ be a right triangle and let $I$ be its incenter (see Fig.\ref{fig:Erone}). Since $\hat B$ is a right angle and since the incircle
is tangent perpendicularly to the three sides of $ABC$, we have  $r=IJ=IK=IH=BJ=BK$, where $r$ is the inradius. The internal bisectors of $ABC$ are
concurrent in $I$ and this implies $AJ=AH=s$ and $CH=CK=t$.

Denoting the area of a triangle with vertical bars we have

$$|ABC|=|AIB|+|BIC|+|CIA|.$$

The half-perimeter of $ABC$ is $p=r+s+t$ while all the three triangles on the r.h.s. of the above equality have altitude $r$
with respect to their sides $AB$, $BC$, and $AC$. Therefore we have

\begin{equation}
\label{AREATRE}
|ABC|=r(r+s+t).
\end{equation}

\begin{rem} In case $ABC$ is not a right triangle, Formula \eqref{AREATRE} generalizes to

\begin{equation}
\label{AREATREGENERAL}
|ABC|=R(r+s+t),
\end{equation}

\noindent where $R$ is the incircle of $ABC$, meaning that $|ABC|$ is a symmetric function in $r,s,t$.
\end{rem}

\begin{teo}\label{teo:right triangle area}
The area of a right triangle $ABC$ is equal to the area of the rectangle of sides $s=AH$ and $t=CH$, where $H$ is the point where the incircle is tangent to the hypotenuse $AC$.

\end{teo}

{\bf Proof} \  Clearly $|ABC|= \displaystyle{{1 \over 2} (s+r)(t+r)}$, and by (\ref{AREATRE}) we get

$$s t +r s +r t + r^2 = 2 (r^2 + r s + r t)$$

\noindent which we simplify into $s t = r^2 + r s + r t = r(r+s+t) = |ABC|$. $\checkmark$

\begin{teo} {\bf [Heron's formula for right triangles]}\label{teo:Heronright} \
The area of a right triangle is equal to  $\sqrt{p(p-a)(p-b)(p-c)}$, where $a$, $b$, $c$ are the
lengths of its sides (taken in any order) and
$p=(a+b+c)/2$  is its half perimeter.

\end{teo}

\smallskip
\noindent
{\bf Proof} \  Using the same notation (as in Fig.1) we need to show that

$$|ABC|^2 = s t r (r+s+t),$$

\noindent  but this is immediate, being $|ABC|=s t$ by Theorem \ref{teo:right triangle area}, and also
$|ABC|=r(r+s+t)$ by (\ref{AREATRE}). $\checkmark$\smallskip

\noindent Theorem \ref{teo:Heronright} shows that, in any right angle triangles, Heron's formula can be rediscovered by starting from the symmetric formula provided by Theorem \ref{teo:right triangle area}. We wish now to extend such a result to any triangle.

\begin{teo} {\bf [Pythagoras]} \
In a right triangle the area of the square whose side is the hypotenuse  is equal to the sum of the areas of
the squares whose sides are the two legs.

\end{teo}

\smallskip
\noindent
{\bf Proof} \  We have $r (r+s +t) = s t$. Multiplying by $2$ and adding $s^2+t^2$ on
both sides we get $s^2+t^2+2 r^2+2 r s +  2 r t  = s^2+t^2 + 2 s t$, i.e.

$$(s+r)^2+(r+t)^2 = (s+t)^2. \quad \checkmark$$

\smallskip

\begin{center}
\begin{figure}[htb]
\centering
\includegraphics[scale=1,viewport=57 600 450 720, clip=true]{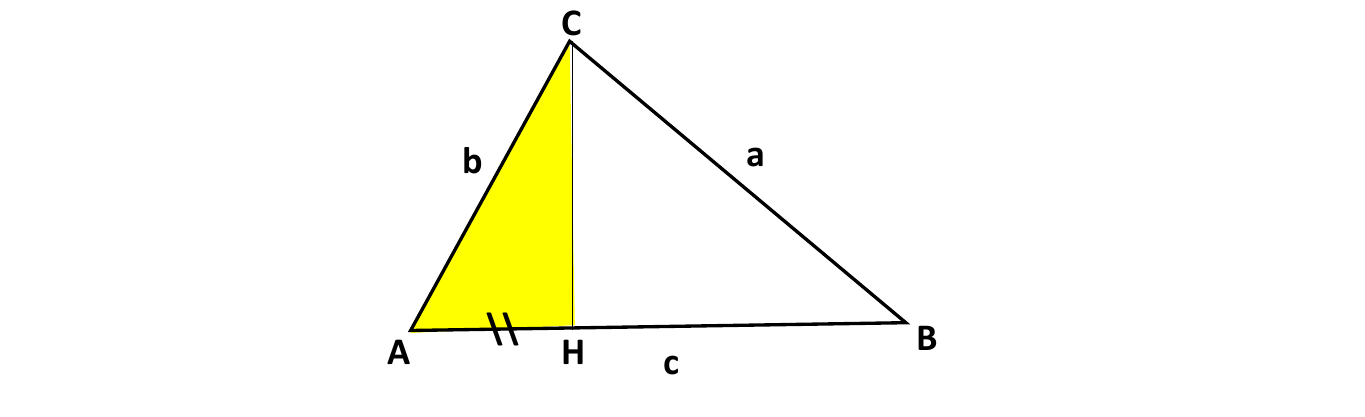}\caption{A generic triangle $ABC$.}\label{fig:EroneGeneric}
\end{figure}
\end{center}

Thanks to Pythagoras' Theorem we can easily extend Heron's Theorem to any triangle. To this, let $ABC$ be a generic triangle, and let $CH$ be the altitude on its edge $AB$ (see Figure \ref{fig:EroneGeneric}), and assume $H$ is between $A$ and $B$ (in any triangle surely exists  an altitude with this property).

Let $p=r+s+t$ be the semiperimeter of $ABC$, and let $\varphi=\sqrt{rst(r+s+t)}=\sqrt{p(p-a)(p-b)(p-c)}$. By Pythagoras’ Theorem  in $AHC$ and $CHB$, we have

$$\begin{array}{l}c^2=(AH+HB)^2=AH^2+BH^2+2(AH)(HB) =\\
= a^2+b^2-2CH^2+2(AH)(HB)=a^2+b^2-2CH^2+2\sqrt{(a^2-CH^2)(b^2-CH^2)}
\end{array}$$

\noindent and, solving for $CH$ we get

$$CH=\frac{\sqrt{4a^2b^2-(c^2-a^2-b^2)^2}}{2c}=\frac{2\sqrt{p(p-a)(p-b)(p-c)}}{c}=\frac{2}{c}\varphi.$$

Therefore, from $2|ABC|=cCH$, Heron’s theorem for $ABC$ follows.\quad \checkmark

\begin{rem} By \eqref{AREATREGENERAL} and Heron's Theorem written in the form $|abc|=\sqrt{rst(r+s+t})$ we can obtain the incircle $R$ of any triangle as a symmetric function of $r,s,t$ as follows (see Figure \ref{fig:HeronGenericTriangle})

\begin{center}
\centering
\begin{figure}[htbp]
\includegraphics[scale=0.68,viewport=-110 580 400 750,clip=true]{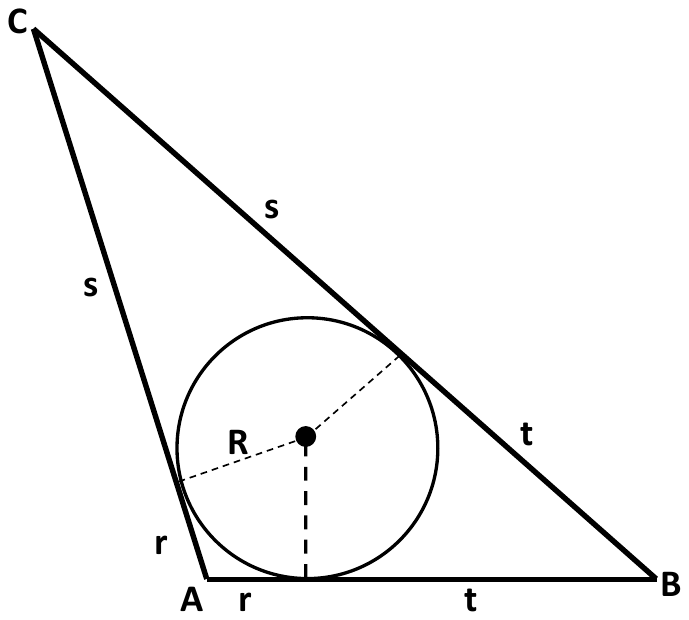}\caption{Incircle of a generic triangle $ABC$.}\label{fig:HeronGenericTriangle}
\end{figure}
\end{center}

\begin{equation}\label{eq:incircle}
R=\sqrt{\frac{rst}{p}}.
\end{equation}

\noindent being $p=r+s+t$ the half perimeter of $ABC$.
\end{rem}

\section{Brahmagupta's formula}

We wish now to show how Brahmagupta's formula can be obtained by exploiting the same idea of symmetry considered in the previous section. with respect to  considered a symmetry , we prove the following result.\smallskip

\begin{teo}\label{teo:areaproduct}
Let $ABC,ADC$ be two triangles inscribed in a same circumference. If $s_1,t_1$ and $s_2,t_2$ are the lengths of the two segments split on the common edge $AC$ by the respective incircles, then

$$|ABC||ACD|=s_1s_2t_1t_2.$$
\end{teo}

\proof In the cyclic quadrangle $ABCD$ the halves of $A\hat{B}C$ and $A\hat{D}C$ are complementary angles. Therefore the shaded right triangles in Figure \ref{fig:Brahmagupta1} are similar, and consequently $\frac{R_1}{r_1}=\frac{r_2}{R_2}$, that is $R_1R_2=r_1r_2$.

\begin{center}
\centering
\begin{figure}[htb]
\includegraphics[scale=0.6,bb=-100 400 480 717,clip=true]{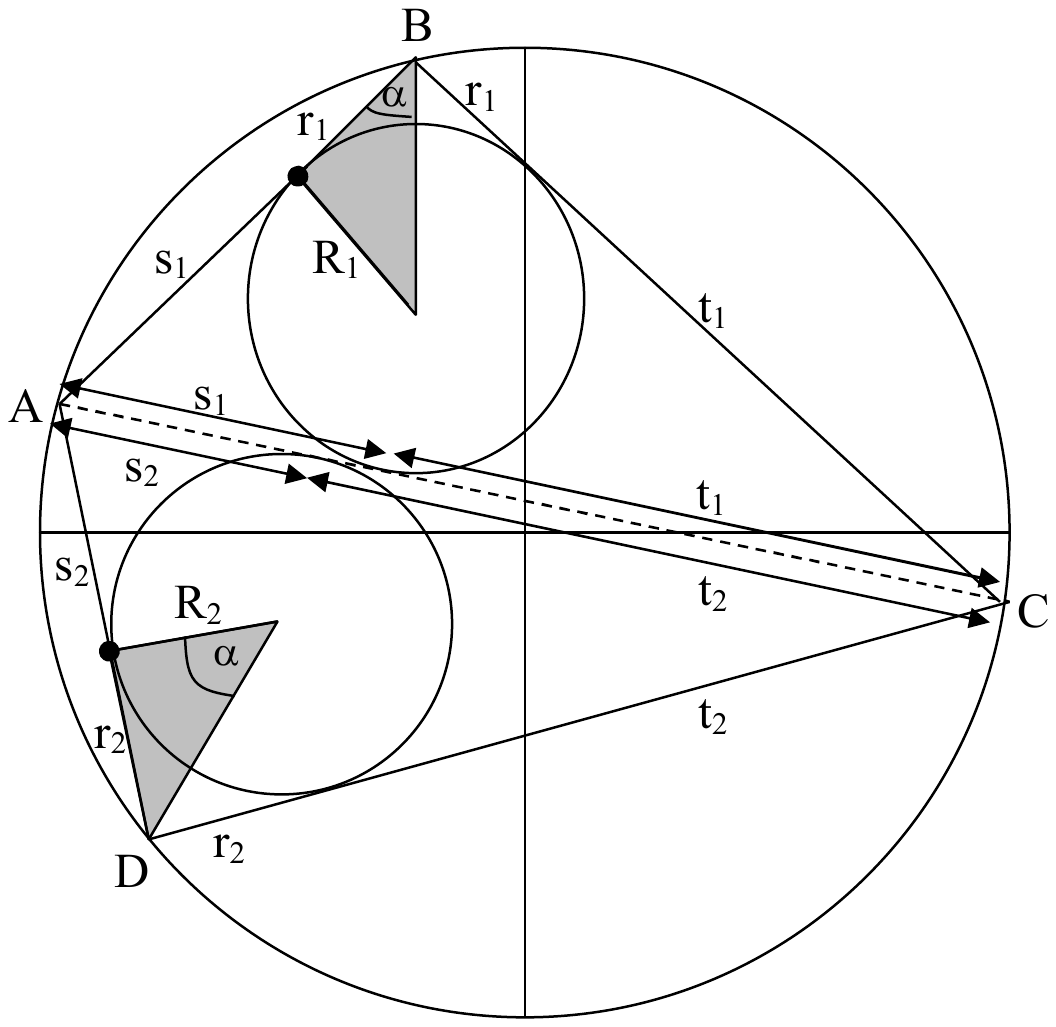}\caption{A general cyclic quadrangle $ABCD$.}\label{fig:Brahmagupta1}
\end{figure}
\end{center}

By \eqref{eq:incircle}  we have $p_1R_1^2=r_1s_1t_1$ and $p_2R_2^2=r_2s_2t_2$, being $p_1,p_2$ the half perimeter of $ABC$ and $ADC$ respectively, so that

$$(p_1R_1p_2R_2)^2=p_1(r_1s_1t_1)p_2(r_2s_2t_2)=p_1p_2(r_1r_2)(s_1t_1s_2t_2)=p_1p_2(R_1R_2)(s_1t_1s_2t_2),$$

\noindent and consequently $|ABC||ACD|=p_1R_1p_2R_2=s_1t_1s_2t_2$. $\checkmark$

\begin{rem}\label{re:areaproduct}\rm Since $R_1R_2=r_1r_2$ we have also $|ABC||ACD|=p_1R_1p_2R_2=p_1p_2r_1r_2$.
\end{rem}\smallskip

\begin{teo}\label{teo:Brahmagupta} {\bf [Brahmagupta's formula for cyclic quadrangles]} \
The area of a convex quadrangle that can be inscribed in a circle
(a cyclic quadrangle) is equal to  $\sqrt{(p-a)(p-b)(p-c)(p-d)}$, where $a$, $b$, $c$, $d$ are the
lengths of its sides (taken in any order) and
$p=(a+b+c+d)/2$  is its half perimeter.
\end{teo}

\proof Let $ABCD$ be split into $ABC$ and $ACD$, as in Figure \ref{fig:Brahmagupta1}, and assume
$a=s_1+r_1=AB,b=t_1+r_1=BC,c=t_2+r_2=CD,d= s_2+r_2 = AD$, so that $p=r_1+r_2+(s_1+t_1)=r_1+r_2+(s_2+t_2)$. Starting from $|ABCD|^2=(|ABC|+|ACD|)^2=|ABC|^2+|ACD|^2+2|ABC||ACD|$, we use the previous theorem, and the Remark, to write $2|ABC||ACD|=s_1t_1s_2t_2+p_1p_2r_1r_2$, where $p_1,p_2$ are the half perimeters of $ABC$ and $ACD$, respectively. Moreover, by Heron's formula, $|ABC|^2=(r_1+s_1+t_1)r_1s_1t_1$, and $|ACD|^2=(r_2+s_2+t_2)r_2s_2t_2$. Then, also using $s_1+t_1=s_2+t_2$, we get

$$\begin{array}{lll}|ABCD|^2&=&\boxed{(r_1+s_1+t_1)r_1s_1t_1}+\boxed{\boxed{(r_2+s_2+t_2)r_2s_2t_2}}+\\
&&+\boxed{\boxed{s_1t_1s_2t_2}}+\boxed{r_1r_2(r_1+s_1+t_1)(r_2+s_2+t_2)}\\\\
&=&\boxed{r_1(r_1+s_1+t_1)(r_2(r_2+s_2+t_2)+s_1t_1)}\\
&&+\boxed{\boxed{s_2t_2(r_2(r_2+s_2+t_2)+s_1t_1)}}\\\\
&=&(r_2(r_2+\boxed{s_2+t_2})+s_1t_1)(r_1(r_1+\boxed{\boxed{s_1+t_1}})+s_2t_2)\\
&=&(r_2(r_2+\boxed{\boxed{s_1+t_1}})+s_1t_1)(r_1(r_1+\boxed{s_2+t_2})+s_2t_2)\\
&=&(r_2+t_1)(r_2+s_1)(r_1+s_2)(r_1+t_2)\\
&=&(p-a)(p-b)(p-c)(p-d).\checkmark
\end{array}$$

\section{The area of a circular polygon having an arbitrary number of edges}

In this section we generalize the previous results to a cyclic polygon $P_n$, having $n+2$ edges for any $n\geq 1$. Let us observe that Heron's formula has been extended to Brahmagupta's formula by considering a cyclic quadrilateral $Q$ as the union of two triangles, $Q=T_1\cup T_2$, and then focusing on the segments $r_1,s_1,t_1$ and $r_2,s_2,t_2$ determined, respectively, on the edges of $T_1$ and $T_2$ by the tangent points of the corresponding incircles. This provides the square of the area of $Q$ as a polynomial function, symmetric under the exchange of $r_1,s_1,t_1$ with $r_2,s_2,t_2$.

As a consequence we are inspired to investigate the square of the area $A(n)$ of a generic cyclic polygon $P_n$ by looking at the partitions of the edges of $P_n$ determined by the tangent points of the incircles of some triangulatio. To this, let us first observe that $P_n$ can always be assumed as a union of $n$ consecutive triangles $T_1,T_2,...,T_n$, all having a common vertex. For $i=1,...,n-1$, denote by $L_{i,i+1}$ the common edge between the two consecutive triangles $T_i,T_{i+1}$. Let $p_j, A_j, R_j$ be, respectively, the semiperimeter, the area and the radius of the incircle of $T_j$, $j=1,...,n$. Also, let $r_j,s_j,t_j$ be the segments cut on the edges of $T_j$ by its incircle, where $L_{i,i+1}=s_i+t_i=s_{i+1}+r_{i+1}$, $i=1,...,n-1$ (see Figure \ref{fig:Catenatriangoli}).

\begin{center}
\centering
\begin{figure}[htb]
\includegraphics[scale=0.688,bb=-50 300 550 800,clip=true]{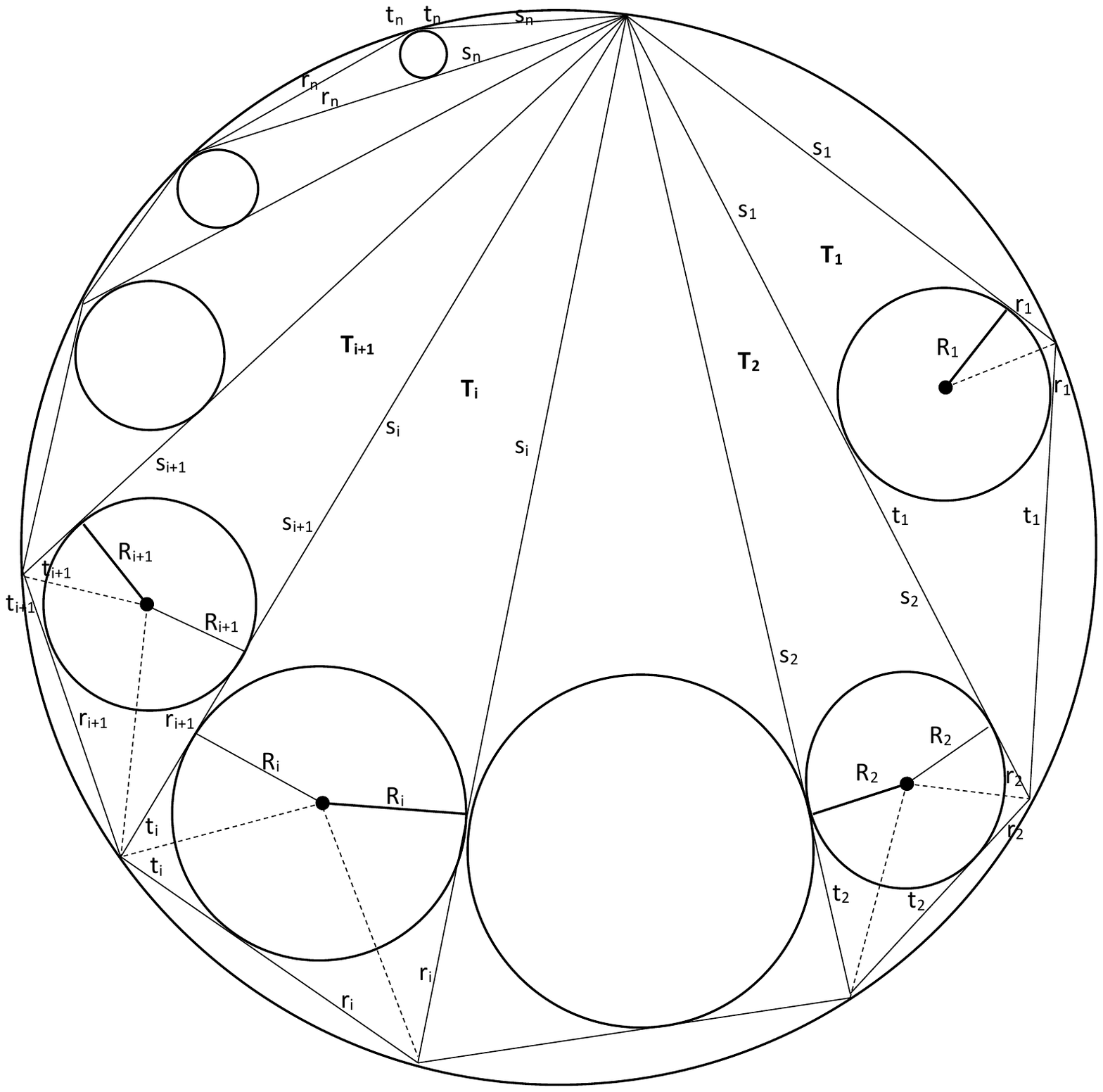}\caption{Triangulation of a cyclic polygons.}\label{fig:Catenatriangoli}
\end{figure}
\end{center}


For the sake of brevity, and in order to avoid heavy notations, in the following theorems we assume to be equal to $1$ all the meaningless products.

\begin{teo}\label{teo:chain}
Let $P_n$ be a cyclic polygon consisting of $n+2$ edges, $n\geq 1$. Then it results


\begin{equation}
\label{eq:chain1} A_hA_k=s_ht_hs_kr_k\prod_{i=h+1}^{k-1}\frac{s_i}{p_i}=p_hr_hp_kt_k\prod_{i=h+1}^{k-1}\frac{p_i}{s_i},\:\:\:\:\:\:\:\:\text{for $1\leq h<k\leq n$}.
\end{equation}

\end{teo}

\proof In order to prove the first equality in \eqref{eq:chain1} we apply Theorem \ref{teo:areaproduct} iteratively, so that

$$\begin{array}{l}A_hA_{h+1}=s_ht_hs_{h+1}r_{h+1}\\
                  A_{h+1}A_{h+2}=s_{h+1}t_{h+1}s_{h+2}r_{h+2}\\
                ...\\
                A_{k-1}A_k=s_{k-1}t_{k-1}s_kr_k.
\end{array}$$

\noindent By multiplying on both sides we get

$$A_h(A_{h+1}A_{h+2}...A_{k-1})^2A_k=s_ht_h\left(\prod_{i=h+1}^{k-1}r_is_i^2t_i\right)s_kr_k.$$

\noindent By Heron's Theorem applied to $T_{h+1}, T_{h+2},...,T_{k-1}$ we have

$$A_hA_k=\frac{s_ht_h\left(\prod_{i=h+1}^{k-1}r_is_i^2t_i\right)s_kr_k}{\prod_{i=h+1}^{k-1}r_is_it_ip_i}=s_ht_hs_kr_k\prod_{i=h+1}^{k-1}\frac{s_i}{p_i},$$

\noindent and the first equality in \eqref{eq:chain1} is obtained. For the proof of the second equality, let us observe that, by similitude, it results $\displaystyle{\frac{R_i}{r_i}=\frac{t_{i+1}}{R_{i+1}}}$, for all $i=1,...,n-1$. Therefore we get

$$\begin{array}{l}R_hR_{h+1}=r_ht_{h+1}\\
                  R_{h+1}R_{h+2}=r_{h+1}t_{h+2}\\
                ...\\
                R_{k-1}R_k=r_{k-1}t_k.
\end{array}$$

\noindent By multiplying on both sides, it results

$$R_h(R_{h+1}R_{h+2}...R_{k-1})^2R_k=R_hR_k\prod_{i=h+1}^{k-1}R_i^2=r_ht_k\prod_{i=h+1}^{k-1}r_it_i,$$

\noindent and applying \eqref{eq:incircle} to all $R_i^2$ we have

$$R_hR_k=\frac{r_ht_k\displaystyle{\prod_{i=h+1}^{k-1}r_it_i}}{\displaystyle{\prod_{i=h+1}^{k-1}\frac{r_is_it_i}{p_i}}},$$

\noindent and consequently

\begin{equation}\label{eq:chain2}
R_hR_k=r_ht_k\prod_{i=h+1}^{k-1}\frac{p_i}{s_i}\:\:\:\:\text{for $1\leq h<k\leq n$}.
\end{equation}

\noindent then The second equality in \eqref{eq:chain1} follows immediately from \eqref{eq:chain2}, being $A_hA_k=p_hR_hp_kR_k$.$\checkmark$\smallskip

Assuming $P=T_1\cup T_2\cup...\cup T_n$ as in Figure \ref{fig:Catenatriangoli}, and using the same notations as above we can now prove a general formula for the area of a cyclic polygon with any number of edges.

\begin{teo}\label{teo:cyclicpolygonarea}
Let $P_n$ be a cyclic polygon consisting of $n+2$ edges, $n\geq 1$, and let $A$ be its area. Then it results

$$A(n)^2=\left(p_1r_1+\sum_{q=2}^nr_qs_q\prod_{m=2}^{q-1}\frac{s_m}{p_m}\right)\left(s_1t_1+\sum_{q=2}^np_qt_q\prod_{m=2}^{q-1}\frac{p_m}{s_m}\right).$$
\end{teo}

\proof Since $P_n$ is the union of $T_1,...,T_n$, by Heron's Theorem, and using both equalities in \eqref{eq:chain1}, we can writewe have

$$\begin{array}{l}A^2=(A_1+A_2+...+A_n)^2=\\
\displaystyle{=\sum_{j=1}^nA_j^2+2\sum_{1\leq h<k\leq n}A_hA_k=}\\
\displaystyle{=\sum_{j=1}^np_jr_js_jt_j+\sum_{1\leq h<k\leq n}s_ht_hs_kr_k\prod_{i=h+1}^{k-1}\frac{s_i}{p_i}+\sum_{1\leq h<k\leq n}p_hr_hp_kt_k\prod_{i=h+1}^{k-1}\frac{p_i}{s_i}}.
\end{array}$$

\noindent Let's collect as follows (where each one of the three terms appearing in each bracket comes from the corresponding sum)

$$\begin{array}{l}
\displaystyle{A^2=p_1r_1\left(s_1t_1+0+\sum_{k>1}p_kt_k\prod_{i=2}^{k-1}\frac{p_i}{s_i}\right)+}\\
\displaystyle{+r_2s_2\left(p_2t_2+s_1t_1+ \frac{p_2}{s_2}\sum_{k>2}p_kt_k\prod_{i=3}^{k-1}\frac{p_i}{s_i}\right)+}\\
\displaystyle{+r_3s_3\frac{s_2}{p_2}\left(\frac{p_2}{s_2}p_3t_3+(s_1t_1+p_2t_2)+\frac{p_2}{s_2}\frac{p_3}{s_3}\sum_{k>3}p_kt_k\prod_{i=4}^{k-1}\frac{p_i}{s_i}\right)+}\\
\displaystyle{+r_4s_4\frac{s_2}{p_2}\frac{s_3}{p_3}\left(\frac{p_2}{s_2}\frac{p_3}{s_3}p_4t_4+(s_1t_1+p_2t_2+\frac{p_2}{s_2}p_3t_3)+\frac{p_2}{s_2}\frac{p_3}{s_3}\frac{p_4}{s_4}\sum_{k>4}p_kt_k\prod_{i=5}^{k-1}\frac{p_i}{s_i}\right)+}\\
+...+\\
\displaystyle{+r_ns_n\frac{s_2}{p_2}\frac{s_3}{p_3}...\frac{s_{n-1}}{p_{n-1}}\left(\frac{p_2}{s_2}\frac{p_3}{s_3}...\frac{p_{n-1}}{s_{n-1}}p_nt_n+\left(s_1t_1+p_2t_2+p_3t_3\frac{p_2}{s_2}+...
+p_{n-1}t_{n-1}\frac{p_2}{s_2}\frac{p_3}{s_3}...\frac{p_{n-2}}{s_{n-2}}\right)+0\right)=}\\\\
\displaystyle{=\left(p_1r_1+\sum_{q=2}^nr_qs_q\prod_{m=2}^{q-1}\frac{s_m}{p_m}\right)\left(s_1t_1+\sum_{q=2}^np_qt_q\prod_{m=2}^{q-1}\frac{p_m}{s_m}\right)}.
\end{array}$$

\hfill$\checkmark$\smallskip

\begin{rem}\label{re:semiperimeter} We emphasize that the formula obtained for $A(n)^2$ is symmetric under mutually exchanging $r_j$ with $t_j$, and $s_j$ with $p_j$, for all $j\in\{1,...,n\}$. We also note that only the terms concerning the partitions of the edges of the polygon $P_n$ are in fact necessary. Indeed, $p_1=r_1+s_1+t_1$ can be immediately computed once we know the terms $r_1, s_1, t_1$ determined on the edges of $P_n$ by the incircle of $T_1$. For $1<q<n-1$, due to consecutiveness, we have $s_q=s_{q-1}+t_{q-1}-r_q$, so that, by recursion, we get

\begin{eqnarray}
\label{eq:s_parts} \displaystyle{s_q=s_1+\sum_{h=1}^{q-1}t_h-\sum_{k=2}^{q}r_k,}\\
\label{eq:semiperimeters} p_q=r_q+s_q+t_q=&\left\{\begin{array}{ll}\displaystyle{s_1+t_1+t_2} &\text{if $q=2$}\\
                                                                 \displaystyle{s_1+\sum_{h=1}^{q}t_h-\sum_{k=2}^{q-1}r_k}&\text{if $2<q<n-1$}.
                                                                 \end{array}\right.
\end{eqnarray}

\noindent Consequently, all terms appearing in $A(n)^2$ can be computed once $s_1, s_n, r_j,t_j$ are known.
\end{rem}

\noindent\textbf{Examples}

\noindent We can easily rediscover Heron's and Brahmagupta's results from the provided symmetric function.

\noindent$\bullet$ $n=1$

$$A(1)^2=\left(p_1r_1+0\right)\left(s_1t_1+0\right)=p_1r_1s_1t_1\:\:\:\:\textbf{Heron's theorem}.$$

\noindent$\bullet$ $n=2$

$$\begin{array}{l}A(2)^2=(p_1r_1+r_2s_2)(s_1t_1+p_2t_2)=\\
=((r_1+s_1+t_1)r_1+r_2s_2)(s_1t_1+(r_2+s_2+t_2)t_2=\\
=((r_1+s_2+r_2)r_1+r_2s_2)(s_1t_1+(s_1+t_1+t_2)t_2=\\
=(r_1+s_2)(r_1+r_2)(t_2+s_1)(t_2+t_1).
\end{array}$$

\noindent Being $s_1+t_1=s_2+r_2$, it is $p=r_1+r_2+s_1+t_1=r_1+r_2+t_2+s_2$, so that, if we denote, respectively, by $a,b,c,d$ the edges $s_2+t_2$, $r_2+t_2$, $r_1+s_1$ and $r_1+t_1$, then

$$A(2)^2=(p-a)(p-b)(p-c)(p-d)\:\:\:\:\textbf{Brahmagupta's theorem}.$$

\noindent$\bullet$ $n=3$ \textbf{generalization of Brahmagupta's theorem to cyclic pentagons}

$$A(3)^2=\left(p_1r_1+r_2s_2+r_3s_3\frac{s_2}{p_2}\right)\left(s_1t_1+p_2t_2+p_3t_3\frac{p_2}{s_2}\right).$$

\noindent$\bullet$ $n=4$ \textbf{generalization of Brahmagupta's theorem to cyclic hexagons}

$$A(4)^2=\left(p_1r_1+r_2s_2+r_3s_3\frac{s_2}{p_2}+r_4s_4\frac{s_2}{p_2}\frac{s_3}{p_3}\right)\left(s_1t_1+p_2t_2+p_3t_3\frac{p_2}{s_2}+p_4t_4\frac{p_2}{s_2}\frac{p_3}{s_3}\right).$$\bigskip

We can even extend the formula to $n=0$ by assuming $A(0)=0$, where $P_0$ is a polygon degenerated in a segment, which can be obtained by progressively removing an edge from a starting polygon $P_n$ having $n+2$ edges

\section{Conclusion and remarks}
We have shown that Heron's and Brahmagupta's theorems can be extended to a formula that provides the square of the area of a convex cyclic polygon as a symmetric polynomial of the segments determined on the edges by the incircles of a suitable triangulation. We remark that the formula is coordinate-free as one should expect from the intrinsic geometric nature of the problem. Differently, using for instance Green's theorem, it would be quite easy to provide a coordinate dependent result.

In our opinion the obtained formula is the natural generalization of what happens for triangle and cyclic quadrilaterals, where the lengths of the edges explicitly appear in the computation of the area. This is just because the number of involved edges is small, so that the segments determined by the edge partitions induced by the incircles can be easily related to the original lengths of the edges of the polygon. We also remark that the incircles can be elementary constructed, so that the provided formula also determines an elementary computation of the square of the area of any convex cyclic polygon

\bibliographystyle{plain}
\bibliography{Heron_Brahmagupta}

\begin{thebibliography}{10}

\bibitem{CG}
H.~S.~M. Coxeter and S.~L. Greitzer.
\newblock {\em Geometry revisited}, volume~19 of {\em New Mathematical
  Library}.
\newblock Random House, Inc., New York, 1967.

\bibitem{CK}
G\'{a}bor Cz\'{e}dli and \'{A}d\'{a}m Kunos.
\newblock Geometric constructibility of cyclic polygons and a limit theorem.
\newblock {\em Acta Sci. Math. (Szeged)}, 81(3-4):643--683, 2015.

\bibitem{D}
Jason DeBlois.
\newblock The geometry of cyclic hyperbolic polygons.
\newblock {\em Rocky Mountain J. Math.}, 46(3):801--862, 2016.

\bibitem{KSS}
Hana Kou\v{r}imsk\'{a}, Lara Skuppin, and Boris Springborn.
\newblock A variational principle for cyclic polygons with prescribed edge
  lengths.
\newblock In {\em Advances in discrete differential geometry}, pages 177--195.
  Springer, [Berlin], 2016.

\bibitem{MRR}
F.~Miller Maley, David~P. Robbins, and Julie Roskies.
\newblock On the areas of cyclic and semicyclic polygons.
\newblock {\em Adv. in Appl. Math.}, 34(4):669--689, 2005.

\bibitem{DG}
Dao~Thanh Oai and Leonard~Mihai Giugiuc.
\newblock The new inequality in a cyclic polygon.
\newblock {\em Int. J. Geom.}, 6(1):5--8, 2017.

\bibitem{PAK}
I.~Pak.
\newblock The area of cyclic polygons: Recent progress on robbins' conjectures.
\newblock {\em Advances in Applied Mathematics}, 34(4):690 -- 696, 2005.
\newblock Special Issue Dedicated to Dr. David P. Robbins.

\bibitem{Robbins1}
D.~P. Robbins.
\newblock Areas of polygons inscribed in a circle.
\newblock {\em Discrete {\&} Computational Geometry}, 12(2):223--236, Dec 1994.

\bibitem{Robbins2}
D.~P. Robbins.
\newblock Areas of polygons inscribed in a circle.
\newblock {\em The American Mathematical Monthly}, 102(6):523--530, 1995.

\bibitem{V}
V.~V. Varfolomeev.
\newblock Inscribed polygons and {H}eron polynomials.
\newblock {\em Mat. Sb.}, 194(3):3--24, 2003.

\bibitem{WX}
Shasha Wang and Wen-Qing Xu.
\newblock Random cyclic polygons from {D}irichlet distributions and
  approximations of {$\pi$}.
\newblock {\em Statist. Probab. Lett.}, 140:84--90, 2018.

\end{thebibliography}


\end{document}